\documentclass[12pt]{article}

\newcommand{\be}{\begin{equation}}
\newcommand{\ee}{\end{equation}}
\newcommand{\bea}{\begin{eqnarray}}
\newcommand{\eea}{\end{eqnarray}}
\newcommand{\barray}{\begin{array}}
\newcommand{\earray}{\end{array}}
\newcommand{\pa}{\partial}
\newcommand{\nn}{\nonumber}
\newcommand{\bitem}{\begin{itemize}}
\newcommand{\eitem}{\end{itemize}}
\newtheorem{teo}{Theorem}[section]
\newcommand{\bt}{\begin{teo}}
\newcommand{\et}{\end{teo}}
\newtheorem{Def}{Definition}[section]
\newcommand{\bd}{\begin{Def}}
\newcommand{\ed}{\end{Def}}
\newtheorem{lem}{Lemma}[section]
\newcommand{\bl}{\begin{lem}}
\newcommand{\el}{\end{lem}}
\newtheorem{prop}{Proposition}[section]
\newcommand{\bp}{\begin{prop}}
\newcommand{\ep}{\end{prop}}
\newtheorem{cor}{Corollary}[section]
\newcommand{\bc}{\begin{cor}}
\newcommand{\ec}{\end{cor}}
\newtheorem{ex}{Example}[section]
\newcommand{\bex}{\begin{ex}}
\newcommand{\eex}{\end{ex}}
\newtheorem{rem}{Remark}[section]
\newcommand{\br}{\begin{rem}}
\newcommand{\er}{\end{rem}}

\catcode `\@=11
\@addtoreset{equation}{section}

\begin{document}

\begin{center}
{\Large \textbf{Nonlocal Hamiltonian operators
of hydrodynamic \\ type with flat metrics,
integrable hierarchies and \\ the equations of
associativity\footnote{This work was supported by the
Alexander von Humboldt Foundation (Germany),
the Russian Foundation for Basic Research
(Grant No. 03-01-00782), and the program of support for
the leading scientific schools (Grant No. 2185.2003.1).}}}
\end{center}

\medskip

\begin{center}
{\large {O. I. Mokhov}}
\end{center}

\bigskip

\section{Introduction} \label{vved}

In this paper we solve the problem of describing all
nonlocal Hamiltonian operators of hydrodynamic type
with flat metrics and establish that this
nontrivial special class of Hamiltonian operators
is closely connected with the associativity equations of
two-dimensional topological quantum field theories
and the theory of Frobenius manifolds.
It is shown that the Hamiltonian operators of this class
are of special interest for many other reasons too.
In particular, we prove in this paper that any such
Hamiltonian operator always defines integrable
structural flows (systems of hydrodynamic type),
always gives a nontrivial pencil of compatible
Hamiltonian operators and generates integrable
hierarchies of hydrodynamic type.
It is proved that the affinors of any such
Hamiltonian operator generate some special integrals
in involution. The nonlinear systems describing
integrals in involution are of independent
great interest. The equations of associativity
of two-dimensional topological quantum field theories
(the Witten--Dijkgraaf--Verlinde--Verlinde and Dubrovin
equations) describe an important special
class of integrals in involution and a special class
of nonlocal Hamiltonian operators of hydrodynamic type
with flat metrics. It is shown that any
$N$-dimensional Frobenius manifold
can be locally presented by a certain special flat
$N$-dimensional submanifold with flat normal bundle
in a $2N$-dimensional pseudo-Euclidean
space and this submanifold is defined uniquely
up to motions. We will devote a separate paper
to the properties
of this construction and to the properties
of this special class of flat submanifolds
with flat normal bundle (we mean the class corresponding to
Frobenius manifolds).

We recall that the general
nonlocal Hamiltonian operators of
hydrodynamic type, namely, the Hamiltonian operators
of the form
\be
P^{ij} = g^{ij}(u(x)) {d \over dx} +
b^{ij}_k (u(x))\, u^k_x + \sum_{n =1}^L
\varepsilon^n (w_n)^i_k (u (x))
u^k_x \left ( {d \over dx} \right )^{-1} \circ
(w_n)^j_s (u (x)) u^s_x, \label{nonle}
\ee
where $\det (g^{ij} (u)) \neq 0,$
$\varepsilon^n = \pm 1,$ $1 \leq n \leq L,$ $u^1, \ldots, u^N$
are local coordinates, $u = (u^1, \ldots, u^N),$
$u^i (x),$ $1 \leq i \leq N,$ are functions (fields)
of one independent variable $x$, the coefficients
$g^{ij} (u),$ $b^{ij}_k (u),$ $(w_n)^i_j (u),$
$1 \leq i, j, k \leq N,$ $1 \leq n \leq L,$ are smooth
functions of local coordinates,
were studied by Ferapontov in the paper \cite{1}
in connection with vital necessities of
the Hamiltonian theory of systems of
hydrodynamic type
(see also \cite{2}, \cite{2a}).

The Hamiltonian operators of the general form
(\ref{nonle}) (local and nonlocal) play a key role
in the Hamiltonian theory of systems of hydrodynamic type.
We recall that an operator $P^{ij}$ is said to be
{\it Hamiltonian} if the operator defines a Poisson bracket
\be
\{ I, J \} = \int {\delta I \over \delta u^i (x)} P^{ij}
{\delta J \over \delta u^j (x)} dx \label{brac}
\ee
for arbitrary functionals $I$ and $J$ on the space of the fields
$u^i (x)$, i.e. the bracket
(\ref{brac}) is skew-symmetric and
satisfies the Jacobi identity.

It is proved in the paper [1] that
the operator (\ref{nonle}) is Hamiltonian if and only if
$g^{ij} (u)$ is a symmetric (pseudo-Riemannian)
contravariant metric and also the coefficients of the
operator satisfy the following relations:

1) $b^{ij}_k (u) = - g^{is} (u) \Gamma ^j_{sk} (u),$
where
$\Gamma^j_{sk} (u)$ is the Riemannian connection
generated by the contravariant metric
$g^{ij} (u)$
(the Levi-Civita connection),

2) $g^{ik} (u) (w_n)^j_k (u) =
g^{jk} (u) (w_n)^i_k (u),$

3) $\nabla_k (w_n)^i_j (u) =
\nabla_j (w_n)^i_k (u),$ where $\nabla_k$
is the operator of covariant differentiation
generated by the Levi-Civita connection
$\Gamma^j_{sk} (u)$ of the metric $g^{ij} (u),$

4) $R^{ij}_{kl} (u) =  \sum_{n =1}^L
\varepsilon^n
\left ( (w_n)^i_l (u) (w_n)^j_k (u)
- (w_n)^j_l (u) (w_n)^i_k (u) \right ),$ where
$$R^{ij}_{kl} (u)= g^{is} (u) R^j_{skl} (u)$$
is the Riemannian curvature tensor of the metric
$g^{ij} (u),$

5) the family of tensors of type (1,1) (i.e., {\it affinors})
$(w_n)^i_j (u),$ $1 \leq n \leq L,$  is commutative:
$[w_n (u), w_m (u)] =0.$

A Hamiltonian operator of the form
(\ref{nonle}) exactly corresponds to an $N$-dimensional
surface with flat normal bundle embedded in
a pseudo-Euclidean space $E^{N + L}$.
Here, the covariant metric $g_{ij} (u),$ for which
$g_{is} (u) g^{sj} (u) = \delta^j_i,$ is the first
fundamental form, and the affinors
$w_n (u),$ $1 \leq n \leq L,$ are the corresponding
Weingarten operators of this embedded surface
($g_{is} (u) (w_n)^s_j (u)$ are the corresponding
second fundamental forms).
Correspondingly, the relations 2)--4)
are the Gauss--Peterson--Codazzi
equations for an $N$-dimensional surface
with flat normal bundle embedded in a pseudo-Euclidean
space $E^{N + L}$ \cite{1}. The relations 5)
are equivalent to the Ricci equations for
this embedded surface.

Taking into account the further applications
to arbitrary Frobenius manifolds, we prefer
to consider the general nonlocal Hamiltonian operators
of hydrodynamic type in the form
\be
P^{ij} = g^{ij}(u(x)) {d \over dx} +
b^{ij}_k (u(x))\, u^k_x + \sum_{m =1}^L \sum_{n =1}^L
\mu^{mn} (w_m)^i_k (u (x))
u^k_x \left ( {d \over dx} \right )^{-1} \circ
(w_n)^j_s (u (x)) u^s_x, \label{nonlm}
\ee
where $\det (g^{ij} (u)) \neq 0,$
$\mu^{mn}$ is an arbitrary nondegenerate symmetric
constant matrix. It is obvious that
considering linear transformations in the
vector space of the affinors $w_n (u),$ $1 \leq n \leq L,$
i.e., changing in (\ref{nonlm})
all the affinors $w_n (u)$ to
$c^l_n \widetilde w_l (u)$, where $c^l_n$ is
an arbitrary nondegenerate constant matrix,
$w_n (u) = c^l_n \widetilde w_l (u),$
any operator of the form (\ref{nonlm}) can be
reduced
to the form (\ref{nonle}) and conversely.
Among all the conditions 1)--5) for the Hamiltonian property
of the operator
(\ref{nonle}), these transformations change only the condition
4) for the Riemannian curvature tensor of the metric.
The condition 4) for the operator
(\ref{nonlm}) takes the form
$$R^{ij}_{kl} (u) =  \sum_{m =1}^L \sum_{n =1}^L
\mu^{mn}
\left ( (w_m)^i_l (u) (w_n)^j_k (u)
- (w_m)^j_l (u) (w_n)^i_k (u) \right ),$$
all the other conditions 1)--3) and 5) for the Hamiltonian property
do not change.

We write all the relations for the coefficients
of the nonlocal Hamiltonian operator
(\ref{nonlm}) in a form convenient for further use.

\bl
The operator (\ref{nonlm}) is Hamiltonian if and only
if its coefficients satisfy the relations
\be
g^{ij} = g^{ji}, \label{01}
\ee
\be
{\pa g^{ij} \over \pa u^k} = b^{ij}_k + b^{ji}_k, \label{02}
\ee
\be
g^{is} b^{jk}_s = g^{js} b^{ik}_s, \label{03}
\ee
\be
g^{is} (w_n)^j_s = g^{js} (w_n)^i_s, \label{04}
\ee
\be
(w_n)^i_s (w_m)^s_j =
(w_m)^i_s (w_n)^s_j, \label{05}
\ee
\be
g^{is} g^{jr} {\pa (w_n)^k_r \over \pa u^s}
- g^{jr} b^{ik}_s (w_n)^s_r  =
g^{js} g^{ir} {\pa (w_n)^k_r \over \pa u^s}
- g^{ir} b^{jk}_s (w_n)^s_r, \label{06}
\ee
\be
g^{is} \left ( {\pa b^{jk}_s \over \pa u^r}
- {\pa b^{jk}_r \over \pa u^s} \right )
+ b^{ij}_s b^{sk}_r - b^{ik}_s b^{sj}_r =
\sum_{m = 1}^L \sum_{n =1}^L \mu^{mn} g^{is}
\left ( (w_m)^j_r (w_n)^k_s -
 (w_m)^j_s (w_n)^k_r \right ). \label{07}
\ee
\el

\section{Pencil of Hamiltonian operators} \label{pencil}

Let us consider the important special case
of the nonlocal Hamiltonian operators of the form
(\ref{nonlm}) when the metric $g^{ij} (u)$ is flat.
We recall that any flat metric uniquely defines
a local Hamiltonian operator of hydrodynamic type
(i.e., a Hamiltonian operator of the form (\ref{nonlm})
with zero affinors) or, in other words,
a Dubrovin--Novikov Hamiltonian operator \cite{2}.
In this paper, we prove that for any flat metric
there is also a remarkable class of nonlocal
Hamiltonian operators of hydrodynamic type
with this flat metric and nontrivial affinors,
and moreover, these Hamiltonian operators
have important applications in the theory of
Frobenius manifolds and integrable hierarchies.

First of all, we note the following important
property of nonlocal Hamiltonian operators
of hydrodynamic type with flat metrics.
We recall that two Hamiltonian operators
are said to be {\it compatible} if any linear combination
of these Hamiltonian operators is also Hamiltonian \cite{5},
i.e., they form a {\it pencil of Hamiltonian operators}.

\bl  \label{fl}
The metric $g^{ij} (u)$ of a Hamiltonian operator
of the form (\ref{nonlm}) is flat if and only if
this operator defines in fact a pencil of
compatible Hamiltonian operators
\bea
P^{ij}_{\lambda_1, \lambda_2} &=&
\lambda_1 \left ( g^{ij}(u(x)) {d \over dx} +
b^{ij}_k (u(x))\, u^k_x \right ) + \nn\\
&+&
\lambda_2 \sum_{m =1}^L \sum_{n =1}^L
\mu^{mn} (w_m)^i_k (u (x))
u^k_x \left ( {d \over dx} \right )^{-1} \circ
(w_n)^j_s (u (x)) u^s_x, \label{nonlma}
\eea
where $\lambda_1$ and $\lambda_2$ are arbitrary constants.
\el

{\it Proof}. Actually, if the operator (\ref{nonlm})
is Hamiltonian, then its coefficients satisfy the relations
(\ref{01})--(\ref{07}). It is obvious that in this case
the relations (\ref{01})--(\ref{06}) for the operator
(\ref{nonlma}) are always satisfied for any constants
$\lambda_1$ and $\lambda_2$, and the relation (\ref{07})
is satisfied for any constants $\lambda_1$ and $\lambda_2$
if and only if the left and the right parts of this relation
are equal to zero identically.

It follows from the relations (\ref{01})--(\ref{03})
for the Hamiltonian operator
(\ref{nonlm}) that the Riemannian curvature tensor
of the metric $g^{ij} (u)$ has the form
\be
R^{ijk}_r (u) = g^{is} (u) R^{jk}_{sr} (u) =
g^{is} (u) \left ( {\pa b^{jk}_s \over \pa u^r}
- {\pa b^{jk}_r \over \pa u^s} \right )
+ b^{ij}_s (u) b^{sk}_r (u) - b^{ik}_s (u) b^{sj}_r (u). \label{kr}
\ee
Consequently, if the metric $g^{ij} (u)$
of a Hamiltonian operator of the form (\ref{nonlm}) is flat,
i.e., $R^{ijk}_r (u) = 0,$ then the relation (\ref{07})
becomes
$$\sum_{m = 1}^L \sum_{n =1}^L \mu^{mn} g^{is}
\left ( (w_m)^j_r (u) (w_n)^k_s (u) -
 (w_m)^j_s (u) (w_n)^k_r (u) \right ) = 0.$$
Thus, the metric $g^{ij} (u)$ of a Hamiltonian operator
of the form (\ref{nonlm}) is flat if and only if
the left and the right parts of the relation
(\ref{07}) for the Hamiltonian operator
(\ref{nonlm}) are equal to zero identically,
and in this case the left
and the right parts of the relation
(\ref{07}) for the operator
(\ref{nonlma}) are also equal to zero identically
for any constants $\lambda_1$ and $\lambda_2$,
i.e., we get a pencil of compatible
Hamiltonian operators (\ref{nonlma}).
Note also that for any pencil of Hamiltonian
operators $P^{ij}_{\lambda_1, \lambda_2}$ (\ref{nonlma})
it follows immediately from the Dubrovin--Novikov
theorem \cite{2} for the local operator that the metric
$g^{ij} (u)$ is flat. Lemma \ref{fl} is proved.

Thus, if the metric
$g^{ij} (u)$ of a Hamiltonian operator
of the form (\ref{nonlm}) is flat, then the operator
\be
P^{ij}_{0, 1} =
 \sum_{m =1}^L \sum_{n =1}^L
\mu^{mn} (w_m)^i_k (u (x))
u^k_x \left ( {d \over dx} \right )^{-1} \circ
(w_n)^j_s (u (x)) u^s_x \label{nonlmb}
\ee
is also a Hamiltonian operator, and moreover,
in this case, this Hamiltonian operator
is always compatible with the local Hamiltonian
operator of hydrodynamic type
(the Dubrovin--Novikov operator)
\be
P^{ij}_{1, 0} =
g^{ij}(u(x)) {d \over dx} +
b^{ij}_k (u(x))\, u^k_x.  \label{nonlmc}
\ee

The compatible Hamiltonian operators
(\ref{nonlmb}) and
(\ref{nonlmc}) always generate the corresponding
integrable bi-Hamiltonian hierarchies. We construct
these integrable hierarchies further in this paper.

\section{Integrability of structural flows} \label{potoki}

We recall that the systems of hydrodynamic type
\be
u^i_{t_n} = (w_n)^i_j (u) u^j_x, \ \ \ \ 1 \leq n \leq L, \label{str}
\ee
are called {\it structural flows} of the nonlocal
Hamiltonian operator of hydrodynamic type (\ref{nonlm})
(see \cite{1}, \cite{4}).

\bl
For any nonlocal Hamiltonian operator
of hydrodynamic type with a flat metric
all the structural flows (\ref{str})
are always commuting integrable bi-Hamiltonian
systems of hydrodynamic type.
\el

{\it Proof}.
As was proved by Maltsev and Novikov in \cite{4}
(see also \cite{1}), the structural flows of any
nonlocal Hamiltonian operator of
hydrodynamic type (\ref{nonlm}) are always Hamiltonian
with respect to this Hamiltonian operator.
Let us consider an arbitrary nonlocal
Hamiltonian operator of hydrodynamic type
(\ref{nonlm}) with a flat metric $g^{ij} (u)$
and
the pencil of compatible Hamiltonian operators
(\ref{nonlma}) corresponding to this Hamiltonian
operator. The corresponding structural flows
must be Hamiltonian with respect to any of the
operators of the Hamiltonian pencil
(\ref{nonlma}) and, consequently,
they are integrable bi-Hamiltonian systems.

\section{Description of nonlocal Hamiltonian
operators \\ of hydrodynamic type with
flat metrics} \label{nonop}

Let us describe all the nonlocal
Hamiltonian operators of hydrodynamic type
with flat metrics.

The form of the Hamiltonian operator (\ref{nonlm})
is invariant with respect to local changes of coordinates,
and also all the coefficients of the operator are
transformed as the corresponding differential-geometric objects.
Since the metric is flat, there exist local coordinates
in which the metric is reduced to a constant form
$\eta^{ij},$ $\eta^{ij}= {\rm const},$ $\det \eta^{ij} \neq 0,$
$\eta^{ij} = \eta^{ji}$.
In these local coordinates all the coefficients of
the Levi-Civita connection are equal to zero, and
the Hamiltonian operator has the form:
\be
\widetilde P^{ij} = \eta^{ij} {d \over dx} +
\sum_{m =1}^L \sum_{n =1}^L
\mu^{mn} (\widetilde w_m)^i_k (u (x))
u^k_x \left ( {d \over dx} \right )^{-1} \circ
(\widetilde w_n)^j_s (u (x)) u^s_x. \label{nonl3}
\ee

\bt  The operator (\ref{nonl3}), where $\eta^{ij}$ and $\mu^{mn}$
are arbitrary nondegenerate symmetric constant matrices,
is Hamiltonian if and only if
there exist functions $\psi_n (u),$ $1 \leq n \leq L,$
such that
\be
(\widetilde w_n)^i_j (u) = \eta^{is}
{\pa^2 \psi_n \over \pa u^s \pa u^j}, \label{n}
\ee
and also the following relations are fulfilled:
\be
\sum_{m=1}^N \sum_{n=1}^N \eta^{mn} {\pa^2 \psi_j \over
\pa u^i \pa u^m} {\pa^2 \psi_k \over \pa u^n \pa u^l}
=
\sum_{m=1}^N \sum_{n=1}^N \eta^{mn} {\pa^2 \psi_k \over
\pa u^i \pa u^m} {\pa^2 \psi_j \over \pa u^n \pa u^l}, \label{nonl4}
\ee
\be
\sum_{m=1}^L \sum_{n=1}^L \mu^{mn}
{\pa^2 \psi_m \over
\pa u^i \pa u^j} {\pa^2 \psi_n \over \pa u^k \pa u^l} =
\sum_{m=1}^L \sum_{n=1}^L \mu^{mn}
{\pa^2 \psi_m \over
\pa u^i \pa u^k} {\pa^2 \psi_n \over \pa u^j \pa u^l}. \label{nonl5}
\ee
\et

{\it Proof.}
The relations (\ref{01})--(\ref{03}) for any operator
of the form (\ref{nonl3}) are automatically fulfilled,
and the relation (\ref{06}) for any operator of the form
(\ref{nonl3}) has the form
\be
{\pa (\widetilde w_n)^k_r \over \pa u^s}
 =
{\pa (\widetilde w_n)^k_s \over \pa u^r}, \label{06a}
\ee
and, consequently, there locally exist functions
$\varphi_n^i (u),$ $1 \leq i \leq N,$ $1 \leq n \leq L,$
such that
\be
(\widetilde w_n)^i_j (u) = {\pa \varphi_n^i \over \pa u^j}.
\ee
The relation (\ref{04}) becomes
\be
\eta^{is} {\pa \varphi_n^j \over \pa u^s} =
\eta^{js} {\pa \varphi_n^i \over \pa u^s}
\ee
or, equivalently,
\be
{\pa \left ( \eta_{is} \varphi_n^s \right ) \over \pa u^j} =
{\pa \left ( \eta_{js} \varphi_n^s \right ) \over \pa u^i}, \label{04a}
\ee
where $\eta_{ij}$ is the matrix that is the inverse of the matrix
$\eta^{ij}$:
$\eta_{is} \eta^{sj} = \delta^j_i.$
It follows from the relation (\ref{04a}) that there locally
exist functions $\psi_n (u),$ $1 \leq n \leq L,$ such that
\be
\eta_{is} \varphi_n^s = {\pa \psi_n \over \pa u^i}.
\ee
Thus,
\be
\varphi_n^i = \eta^{is} {\pa \psi_n \over \pa u^s},
\ \ \ (\widetilde w_n)^i_j (u) = \eta^{is}
{\pa^2 \psi_n \over \pa u^s \pa u^j}.
\ee
In this case the relations (\ref{05}) and (\ref{07})
become (\ref{nonl4}) and (\ref{nonl5}) respectively.
Theorem is proved.

\bc
If the functions $\psi_n (u),$ $1 \leq n \leq L,$
are a solution of the system of the equations
(\ref{nonl4}), (\ref{nonl5}), then the systems
of hydrodynamic type
\be
u^i_{t_n} = \eta^{is}
{\pa^2 \psi_n \over \pa u^s \pa u^j} u^j_x,
\ \ \ \ 1 \leq n \leq L, \label{str2}
\ee
are always commuting integrable bi-Hamiltonian systems
of hydrodynamic type. Moreover, in this case
the operator
\be
M^{ij}_1 =
 \sum_{m =1}^L \sum_{n =1}^L
\mu^{mn} \eta^{ip} \eta^{jr} {\pa^2 \psi_m \over \pa u^p \pa u^k}
u^k_x \left ( {d \over dx} \right )^{-1} \circ
{\pa^2 \psi_n \over \pa u^r \pa u^s} u^s_x \label{nonlmd}
\ee
is also a Hamiltonian operator,
and also this Hamiltonian operator is always compatible
with the constant Hamiltonian operator
\be
M^{ij}_2 =
\eta^{ij} {d \over dx}.  \label{nonlme}
\ee
\ec

In arbitrary local coordinates we obtain the following
description of all nonlocal Hamiltonian
operators of hydrodynamic type with flat metrics.

\bt  The operator (\ref{nonlm}) with a flat metric $g^{ij} (u)$
is Hamiltonian if and only if
$b^{ij}_k (u) = - g^{is} (u) \Gamma^j_{sk} (u),$
where
$\Gamma^j_{sk} (u)$ is the flat connection
generated by the flat metric $g^{ij} (u)$, and there locally
exist functions $\Phi_n (u),$ $1 \leq n \leq L,$ such that
\be
(w_n)^i_j (u) = \nabla^i \nabla_j \Phi_n, \label{n2}
\ee
and also the following relations are fulfilled:
\be
 \sum_{n=1}^N  \nabla^n \nabla_i \Phi_j
\nabla_n \nabla_l \Phi_k
=  \sum_{n=1}^N  \nabla^n \nabla_i \Phi_k
\nabla_n \nabla_l \Phi_j, \label{nonl4n}
\ee
\be
\sum_{m=1}^L \sum_{n=1}^L \mu^{mn}
\nabla_i \nabla_j \Phi_m \nabla_k \nabla_l \Phi_n =
\sum_{m=1}^L \sum_{n=1}^L \mu^{mn}
\nabla_i \nabla_k \Phi_m \nabla_j \nabla_l \Phi_n , \label{nonl5n}
\ee
where $\nabla_i$ is the operator of
covariant differentiation defined by the flat connection
$\Gamma^i_{jk} (u)$ generated by the metric
$g^{ij} (u)$, $\nabla^i = g^{is} (u) \nabla_s.$
In particular, in this case the operator
\bea
M^{ij}_{\lambda_1, \lambda_2} &=&
\lambda_1 \left ( g^{ij}(u(x)) {d \over dx} -
g^{is} (u (x)) \Gamma^j_{sk} (u(x))\, u^k_x \right ) + \nn\\
&+&
\lambda_2
 \sum_{m =1}^L \sum_{n =1}^L
\mu^{mn} \nabla^i \nabla_k \Phi_m
u^k_x \left ( {d \over dx} \right )^{-1} \circ
\nabla^j \nabla_s \Phi_n u^s_x \label{nonlmdn}
\eea
is a Hamiltonian operator for any constants
$\lambda_1$ and $\lambda_2$, and the systems
of hydrodynamic type
\be
u^i_{t_n} = \nabla^i \nabla_j
\Phi_n u^j_x,
\ \ \ \ 1 \leq n \leq L, \label{str2n}
\ee
are always commuting integrable bi-Hamiltonian systems
of hydrodynamic type.
\et

\section{Integrable hierarchies} \label{hier5}

Let us consider the recursion operator $R^i_j$
corresponding to the compatible Hamiltonian operators
(\ref{nonlmd}) and (\ref{nonlme}):
\be
R^i_j = \left ( M_1 (M_2)^{- 1} \right )^i_j =
\sum_{m =1}^L \sum_{n =1}^L
\mu^{mn} \eta^{ip} {\pa^2 \psi_m \over \pa u^p \pa u^k}
u^k_x \left ( {d \over dx} \right )^{-1} \circ
{\pa^2 \psi_n \over \pa u^j \pa u^s} u^s_x
\left ( {d \over dx} \right )^{-1}. \label{rek}
\ee
Let us apply the constructed recursion operator
(\ref{rek}) to the system of translations
with respect to $x$, i.e., to the system
\be
u^i_t = u^i_x.
\ee
Any system in the hierarchy
\be
u^i_{t_s} = (R^s)^i_j u^j_x, \ \ \ \ s \in {\bf Z},  \label{ierar}
\ee
is a multi-Hamiltonian integrable system of
hydrodynamic type. In particular,
any system of the form
\be
u^i_{t_1} = R^i_j u^j_x,
\ee
i.e., the system
\be
u^i_{t_1} =
\sum_{m =1}^L \sum_{n =1}^L
\mu^{mn} \eta^{ip} {\pa^2 \psi_m \over \pa u^p \pa u^k}
u^k_x \left ( {d \over dx} \right )^{-1} \circ
{\pa^2 \psi_n \over \pa u^j \pa u^s} u^j u^s_x, \label{s1}
\ee
is integrable.

Since
\be
{\pa \over \pa u^r} \left ( {\pa^2 \psi_n \over
\pa u^j \pa u^s} u^j \right ) =
{\pa^3 \psi_n \over \pa u^j \pa u^s \pa u^r} u^j
+ {\pa^2 \psi_n \over \pa u^r \pa u^s} =
{\pa \over \pa u^s} \left ( {\pa^2 \psi_n \over
\pa u^j \pa u^r} u^j \right ),
\ee
there locally exist functions $F_n (u),$ $1 \leq n \leq L,$
such that
\be
{\pa^2 \psi_n \over
\pa u^j \pa u^s} u^j = {\pa F_n \over \pa u^s},
\ \ \ \ F_n = {\pa \psi_n \over \pa u^j} u^j - \psi_n.
\ee
Thus, the system of hydrodynamic type
(\ref{s1}) has a local form:
\be
u^i_{t_1} =
\sum_{m =1}^L \sum_{n =1}^L
\mu^{mn} \eta^{ip} F_n (u)
{\pa^2 \psi_m \over \pa u^p \pa u^k}
u^k_x. \label{s2}
\ee
This system of hydrodynamic type is
bi-Hamiltonian with respect to the compatible
Hamiltonian operators
(\ref{nonlmd}) and (\ref{nonlme}):
\be
u^i_{t_1} =
\sum_{m =1}^L \sum_{n =1}^L
\mu^{mn} \eta^{ip} \eta^{jr} {\pa^2 \psi_m \over \pa u^p \pa u^k}
u^k_x \left ( {d \over dx} \right )^{-1} \left (
{\pa^2 \psi_n \over \pa u^r \pa u^s} u^s_x {\delta H_1
\over \delta u^j (x)} \right ), \label{s3}
\ee
\be
H_1 = \int h_1 (u(x)) dx, \ \ \ \ h_1 (u(x)) =
{1 \over 2} \eta_{ij} u^i (x) u^j (x),
\ee
\be
u^i_{t_1} = \eta^{ij}  {d \over dx}  {\delta H_2
\over \delta u^j (x)}, \ \ \ \ H_2 = \int h_2 (u(x)) dx, \label{s4}
\ee
since in our case there always exists locally a function
$h_2 (u)$ such that
\be
\sum_{m =1}^L \sum_{n =1}^L
\mu^{mn} {\pa^2 \psi_m \over \pa u^j \pa u^k} F_n (u) =
{\pa^2 h_2 \over \pa u^j \pa u^k}.
\ee

Actually, we have
\bea
&&
{\pa \over \pa u^i} \left ( \sum_{m =1}^L \sum_{n =1}^L
\mu^{mn} {\pa^2 \psi_m \over \pa u^j \pa u^k} F_n (u) \right ) =
\sum_{m =1}^L \sum_{n =1}^L
\mu^{mn} {\pa^3 \psi_m \over \pa u^i \pa u^j \pa u^k} F_n (u) + \nn\\
&&
+ \sum_{m =1}^L \sum_{n =1}^L
\mu^{mn} {\pa^2 \psi_m \over \pa u^j \pa u^k} {\pa F_n  \over
\pa u^i} = \sum_{m =1}^L \sum_{n =1}^L
\mu^{mn} {\pa^3 \psi_m \over \pa u^i \pa u^j \pa u^k} F_n (u)
+ \nn\\
&&
+ \sum_{m =1}^L \sum_{n =1}^L
\mu^{mn} {\pa^2 \psi_m \over \pa u^j \pa u^k} {\pa^2 \psi_n  \over
\pa u^i \pa u^s} u^s = \sum_{m =1}^L \sum_{n =1}^L
\mu^{mn} {\pa^3 \psi_m \over \pa u^i \pa u^j \pa u^k} F_n (u)
+ \nn\\
&&
+ \sum_{m =1}^L \sum_{n =1}^L
\mu^{mn} {\pa^2 \psi_m \over \pa u^j \pa u^i} {\pa^2 \psi_n  \over
\pa u^k \pa u^s} u^s, \label{k1}
\eea
where we have used the relation (\ref{nonl5}).
Consequently, by virtue of symmetry
with respect to the indices $i$ and $j$ we obtain
\be
{\pa \over \pa u^i} \left ( \sum_{m =1}^L \sum_{n =1}^L
\mu^{mn} {\pa^2 \psi_m \over \pa u^j \pa u^k} F_n (u) \right ) =
{\pa \over \pa u^j} \left ( \sum_{m =1}^L \sum_{n =1}^L
\mu^{mn} {\pa^2 \psi_m \over \pa u^i \pa u^k} F_n (u) \right ),
\ee
i.e., there locally exist functions $a_k (u),$ $1 \leq k \leq N,$
such that
\be
 \sum_{m =1}^L \sum_{n =1}^L
\mu^{mn} {\pa^2 \psi_m \over \pa u^j \pa u^k} F_n (u) =
{\pa a_k \over \pa u^j}.
\ee
By virtue of symmetry with respect to
the indices $j$ and $k$ we get
\be
{\pa a_k \over \pa u^j} = {\pa a_j \over \pa u^k},
\ee
i.e., there locally exists a function
$h_2 (u)$ such that
\be
a_k (u) = {\pa h_2 \over \pa u^k}.
\ee
Thus,
\be
 \sum_{m =1}^L \sum_{n =1}^L
\mu^{mn} {\pa^2 \psi_m \over \pa u^j \pa u^k} F_n (u) =
{\pa a_k \over \pa u^j} = {\pa^2 h_2 \over \pa u^j \pa u^k}. \label{k11}
\ee

Consider the following equation of the integrable hierarchy
(\ref{ierar}).
\bea
&&
u^i_{t_2} = (R^2)^i_j u^j_x = \nn\\
&&
= \sum_{m =1}^L \sum_{n =1}^L
\mu^{mn} \eta^{ip} {\pa^2 \psi_m \over \pa u^p \pa u^k}
u^k_x \left ( {d \over dx} \right )^{-1} \circ
{\pa^2 \psi_n \over \pa u^j \pa u^s} u^s_x
\left ( {d \over dx} \right )^{-1} \circ
\eta^{jr}  {d \over dx}  {\delta H_2
\over \delta u^r (x)} = \nn\\
&&
= \sum_{m =1}^L \sum_{n =1}^L
\mu^{mn} \eta^{ip} {\pa^2 \psi_m \over \pa u^p \pa u^k}
u^k_x \left ( {d \over dx} \right )^{-1} \circ
{\pa^2 \psi_n \over \pa u^j \pa u^s} u^s_x
\eta^{jr}  {\pa h_2
\over \pa u^r}.
 \label{ierar2}
\eea

We prove that in our case there always exist locally
functions $G_n (u),$ $1 \leq n \leq L,$ such that
\be
{\pa^2 \psi_n \over
\pa u^j \pa u^s} \eta^{jr}  {\pa h_2
\over \pa u^r} = {\pa G_n \over \pa u^s}. \label{gn}
\ee

Actually, we have
\bea
&&
{\pa \over \pa u^p} \left ( {\pa^2 \psi_n \over
\pa u^j \pa u^s} \eta^{jr}  {\pa h_2
\over \pa u^r} \right ) = {\pa^3 \psi_n \over
\pa u^j \pa u^s \pa u^p} \eta^{jr}  {\pa h_2
\over \pa u^r}
+ {\pa^2 \psi_n \over
\pa u^j \pa u^s} \eta^{jr}  {\pa^2 h_2
\over \pa u^r \pa u^p} = \nn\\
&&
= {\pa^3 \psi_n \over
\pa u^j \pa u^s \pa u^p} \eta^{jr}  {\pa h_2
\over \pa u^r}
+ {\pa^2 \psi_n \over
\pa u^j \pa u^s} \eta^{jr} \left (\sum_{k =1}^L \sum_{l =1}^L
\mu^{kl} {\pa^2 \psi_k \over \pa u^r \pa u^p} F_l (u) \right )
= \nn\\
&&
= {\pa^3 \psi_n \over
\pa u^j \pa u^s \pa u^p} \eta^{jr}  {\pa h_2
\over \pa u^r}
+ \sum_{k =1}^L \sum_{l =1}^L
\mu^{kl} \eta^{jr} {\pa^2 \psi_n \over
\pa u^s \pa u^j} {\pa^2 \psi_k \over \pa u^r \pa u^p} F_l (u)
= \nn\\
&&
= {\pa^3 \psi_n \over
\pa u^j \pa u^s \pa u^p} \eta^{jr}  {\pa h_2
\over \pa u^r}
+ \sum_{k =1}^L \sum_{l =1}^L
\mu^{kl} \eta^{jr} {\pa^2 \psi_k \over
\pa u^s \pa u^j} {\pa^2 \psi_n \over \pa u^r \pa u^p} F_l (u)
= \nn\\
&&
= {\pa^3 \psi_n \over
\pa u^j \pa u^s \pa u^p} \eta^{jr}  {\pa h_2
\over \pa u^r}
+ \sum_{k =1}^L \sum_{l =1}^L
\mu^{kl} \eta^{jr} {\pa^2 \psi_k \over
\pa u^s \pa u^r} {\pa^2 \psi_n \over \pa u^j \pa u^p} F_l (u), \label{k2}
\eea
where we have used the relation (\ref{nonl4}) and
the symmetry of the matrix
$\eta^{jr}$. Therefore, it is proved that
the expression under consideration is
symmetric with respect to the indices
$p$ and $s$, i.e.,
\be
{\pa \over \pa u^p} \left ( {\pa^2 \psi_n \over
\pa u^j \pa u^s} \eta^{jr}  {\pa h_2
\over \pa u^r} \right ) =
{\pa \over \pa u^s} \left ( {\pa^2 \psi_n \over
\pa u^j \pa u^p} \eta^{jr}  {\pa h_2
\over \pa u^r} \right ). \label{k22}
\ee
Consequently, there locally exist functions
$G_n (u),$ $1 \leq n \leq L,$ such that
the relation
(\ref{gn}) is fulfilled, and therefore, it is proved
that the second flow of the integrable hierarchy (\ref{ierar})
has the form of a local system of hydrodynamic type
\be
u^i_{t_2} = \sum_{m =1}^L \sum_{n =1}^L
\mu^{mn} \eta^{ip} G_n (u) {\pa^2 \psi_m \over \pa u^p \pa u^k}
u^k_x.
 \label{ierar2l}
\ee

Repeating the foregoing arguments exactly, we prove
by induction that if the functions
$\psi_n (u),$ $1 \leq n \leq L,$ are a solution of
the system of equations
(\ref{nonl4}), (\ref{nonl5}), then for any $s \geq 1$
and for the corresponding function $h_s (u(x))$
(starting from the function
$h_1 (u(x)) = {1 \over 2} \eta_{ij} u^i (x) u^j (x)$)
there always exist locally functions
$F_n^{(s)} (u),$ $1 \leq n \leq L,$ such that
\be
{\pa^2 \psi_n \over
\pa u^j \pa u^p} \eta^{jr}  {\pa h_s
\over \pa u^r} = {\pa F_n^{(s)} \over \pa u^p}, \label{fsn}
\ee
and there always exists locally a function $h_{s+1} (u(x))$
such that
\be
\sum_{m =1}^L \sum_{n =1}^L
\mu^{mn} {\pa^2 \psi_m \over \pa u^j \pa u^k} F_n^{(s)} (u) =
{\pa^2 h_{s+1} \over \pa u^j \pa u^k}.   \label{hsn}
\ee
Above we have proved that this statement
is true for $s = 1$
(in this case, in particular, $F_n^{(1)} = F_n,$ $F_n^{(2)} = G_n$).
In just the same way it is proved that
if this statement is true for $s = K \geq 1$,
then it is true also for $s = K + 1$
(see (\ref{k1})--(\ref{k11}) and (\ref{k2}), (\ref{k22})).
Thus, it is proved that for any $s \geq 1$
the corresponding flow of the integrable hierarchy (\ref{ierar})
has the form of a local system of hydrodynamic type
\be
u^i_{t_s} = \sum_{m =1}^L \sum_{n =1}^L
\mu^{mn} \eta^{ip} F_n^{(s)} (u) {\pa^2 \psi_m \over \pa u^p \pa u^k}
u^k_x.
 \label{ierarsl}
\ee
All the flows of the hierarchy (\ref{ierar}) are commuting
integrable bi-Hamiltonian systems of hydrodynamic type
with an infinite family of local integrals in involution
with respect to both the Poisson brackets:
\be
u^i_{t_s} =
M^{ij}_1 {\delta H_s
\over \delta u^j (x)} = \{ u^i, H_s \}_1, \ \ \ \
H_s = \int h_s (u(x)) dx,
\ee
\be
u^i_{t_s} = M^{ij}_2  {\delta H_{s+1}
\over \delta u^j (x)} = \{ u^i, H_{s+1} \}_2, \ \ \ \
H_{s+1} = \int h_{s+1} (u(x)) dx,
\ee
\be
 \{ H_p, H_r \}_1 = 0, \ \ \ \
\{ H_p, H_r \}_2 = 0,
\ee
and the densities of the Hamiltonians $h_s (u(x))$
are related by the reccurences (\ref{fsn}), (\ref{hsn}),
which are always resolvable in our case.

\section{Locality and integrability of Hamiltonian \\
systems with nonlocal Poisson brackets} \label{lok}

Let the functions
$\psi_n (u),$ $1 \leq n \leq L,$ be a solution
of the system of equations
(\ref{nonl4}), (\ref{nonl5}), i.e., in particular,
the nonlocal operator $M^{ij}_1$ (\ref{nonlmd})
is Hamiltonian. Let us consider the Hamiltonian system
\be
u^i_t =
M^{ij}_1 {\delta H
\over \delta u^j (x)} \label{arb}
\ee
with an arbitrary Hamiltonian of hydrodynamic type
\be
H = \int h (u(x)) dx. \label{arb2}
\ee
As was proved by Ferapontov in \cite{1},
a Hamiltonian system with a nonlocal Hamiltonian
operator of hydrodynamic type
(\ref{nonle}) and a Hamiltonian of hydrodynamic type
is local if and only if the Hamiltonian
is an integral of all the structural flows of the
Hamiltonian operator.
This statement is also true for Hamiltonian operators
of the form (\ref{nonlmb}), and, moreover,
generally for any weakly nonlocal Hamiltonian
operators (see \cite{4}). We prove that
for the nonlocal Hamiltonian operators
$M^{ij}_1$ (\ref{nonlmd}) constructed in this paper
this condition for Hamiltonians is sufficient
for integrability, i.e., all local
Hamiltonian systems (\ref{arb}), (\ref{arb2}) are
integrable bi-Hamiltonian systems.

\bl
The system (\ref{arb}), (\ref{arb2}) is local
if and only if the density $h (u(x))$ of the Hamiltonian
satisfies the equations
\be
{\pa^2 \psi_n \over
\pa u^j \pa u^s} \eta^{jr}  {\pa^2 h
\over \pa u^r \pa u^p}  =
 {\pa^2 \psi_n \over
\pa u^j \pa u^p} \eta^{jr}  {\pa^2 h
\over \pa u^r \pa u^s}. \label{locint}
\ee
\el

{\it Proof.}
Consider a Hamiltonian system (\ref{arb}), (\ref{arb2}):
\be
u^i_t = M^{ij}_1 {\delta H
\over \delta u^j (x)} =
\sum_{m =1}^L \sum_{n =1}^L
\mu^{mn} \eta^{ip} {\pa^2 \psi_m \over \pa u^p \pa u^k}
u^k_x \left ( {d \over dx} \right )^{-1} \left ( \eta^{jr}
{\pa^2 \psi_n \over \pa u^r \pa u^s} u^s_x {\pa h
\over \pa u^j} \right ). \label{hamsyst}
\ee

The system (\ref{hamsyst}) is local if and only if
there locally exist functions $P_n (u),$ $1 \leq n \leq L,$
such that
\be
{\pa^2 \psi_n \over
\pa u^j \pa u^s} \eta^{jr}  {\pa h
\over \pa u^r} = {\pa P_n \over \pa u^s}, \label{pn}
\ee
i.e., if and only if the consistency relation
\be
{\pa \over \pa u^p} \left ( {\pa^2 \psi_n \over
\pa u^j \pa u^s} \eta^{jr}  {\pa h
\over \pa u^r} \right ) =
{\pa \over \pa u^s} \left ( {\pa^2 \psi_n \over
\pa u^j \pa u^p} \eta^{jr}  {\pa h
\over \pa u^r} \right ) \label{sovm}
\ee
is fulfilled.
Then the system (\ref{hamsyst}) takes a local form
\be
u^i_t = M^{ij}_1 {\delta H
\over \delta u^j (x)} =
\sum_{m =1}^L \sum_{n =1}^L
\mu^{mn} \eta^{ip} {\pa^2 \psi_m \over \pa u^p \pa u^k} P_n (u)
u^k_x. \label{hamsyst2}
\ee

The consistency relation (\ref{sovm}) is equivalent to
the equations (\ref{locint}).

\bt
If a Hamiltonian system (\ref{arb}), (\ref{arb2})
is local, i.e., the density $h(u(x))$ of the Hamiltonian
satisfies the equations
(\ref{locint}), then this system is integrable.
\et

{\it Proof.}
In this case the system (\ref{arb}),
(\ref{arb2}) takes the form  (\ref{hamsyst2}),
(\ref{pn}).
Let us prove that there always exists locally a function
$f (u)$ such that
\be
\sum_{m =1}^L \sum_{n =1}^L
\mu^{mn} {\pa^2 \psi_m \over \pa u^j \pa u^k} P_n (u) =
{\pa^2 f \over \pa u^j \pa u^k}.
\ee

Actually, we have
\bea
&&
{\pa \over \pa u^i} \left ( \sum_{m =1}^L \sum_{n =1}^L
\mu^{mn} {\pa^2 \psi_m \over \pa u^j \pa u^k} P_n (u) \right ) =
\sum_{m =1}^L \sum_{n =1}^L
\mu^{mn} {\pa^3 \psi_m \over \pa u^i \pa u^j \pa u^k} P_n (u) + \nn\\
&&
+ \sum_{m =1}^L \sum_{n =1}^L
\mu^{mn} {\pa^2 \psi_m \over \pa u^j \pa u^k} {\pa P_n  \over
\pa u^i} = \sum_{m =1}^L \sum_{n =1}^L
\mu^{mn} {\pa^3 \psi_m \over \pa u^i \pa u^j \pa u^k} P_n (u)
+ \nn\\
&&
+ \sum_{m =1}^L \sum_{n =1}^L
\mu^{mn} {\pa^2 \psi_m \over \pa u^j \pa u^k}
{\pa^2 \psi_n \over
\pa u^i \pa u^p} \eta^{pr}  {\pa h
\over \pa u^r} = \sum_{m =1}^L \sum_{n =1}^L
\mu^{mn} {\pa^3 \psi_m \over \pa u^i \pa u^j \pa u^k} P_n (u)
+ \nn\\
&&
+ \sum_{m =1}^L \sum_{n =1}^L
\mu^{mn} {\pa^2 \psi_m \over \pa u^j \pa u^i} {\pa^2 \psi_n  \over
\pa u^k \pa u^p} \eta^{pr}  {\pa h
\over \pa u^r},
\eea
where we have used the relation (\ref{nonl5}).
Consequently, by virtue of symmetry with
respect to the indices $i$ and $j$ we get
\be
{\pa \over \pa u^i} \left ( \sum_{m =1}^L \sum_{n =1}^L
\mu^{mn} {\pa^2 \psi_m \over \pa u^j \pa u^k} P_n (u) \right ) =
{\pa \over \pa u^j} \left ( \sum_{m =1}^L \sum_{n =1}^L
\mu^{mn} {\pa^2 \psi_m \over \pa u^i \pa u^k} P_n (u) \right ),
\ee
i.e., there locally exist functions
$b_k (u),$ $1 \leq k \leq N,$ such that
\be
 \sum_{m =1}^L \sum_{n =1}^L
\mu^{mn} {\pa^2 \psi_m \over \pa u^j \pa u^k} P_n (u) =
{\pa b_k \over \pa u^j}.
\ee
By virtue of symmetry with respect to the indices
$j$ and $k$ we get
\be
{\pa b_k \over \pa u^j} = {\pa b_j \over \pa u^k},
\ee
i.e., there locally exists a function $f (u)$
such that
\be
b_k (u) = {\pa f \over \pa u^k}.
\ee
Thus,
\be
 \sum_{m =1}^L \sum_{n =1}^L
\mu^{mn} {\pa^2 \psi_m \over \pa u^j \pa u^k} P_n (u) =
{\pa b_k \over \pa u^j} = {\pa^2 f \over \pa u^j \pa u^k}.
\ee
Consequently, the system (\ref{arb}), (\ref{arb2})
in the case under consideration can be presented in the form
\be
u^i_t = \eta^{ij} {\pa^2 f \over \pa u^j \pa u^k} u^k_x =
M^{ij}_2 {\delta F
\over \delta u^j (x)}, \ \ \ \ F = \int f (u) dx, \label{arb3}
\ee
i.e., it is an integrable bi-Hamiltonian system with
the compatible Hamiltonian operators $M^{ij}_1$
(\ref{nonlmd}) and $M^{ij}_2$ (\ref{nonlme}).

\section{Systems of integrals in involution} \label{invol}

The nonlinear equations of the form
(\ref{nonl4}) and (\ref{locint})
are of independent interest, they have
an important significance and
very natural interpretation.

\bl
The nonlinear equations (\ref{nonl4}) are equivalent
to the condition that the integrals
\be
\Psi_n = \int \psi_n (u(x)) dx, \ \ 1 \leq n \leq L,
\ee
are in involution with respect to the Poisson bracket
given by the constant Hamiltonian operator
$M^{ij}_2$ (\ref{nonlme}):
\be
\{ \Psi_n, \Psi_m \} = 0.
\ee
\el

{\it Proof.}
Actually, we have
\be
\{ \Psi_n, \Psi_m \} = \int {\pa \psi_n \over \pa u^i}
\eta^{ij} {d \over d x} {\pa \psi_m \over \pa u^j} dx =
\int {\pa \psi_n \over \pa u^i}
\eta^{ij} {\pa^2 \psi_m \over \pa u^j \pa u^k} u^k_x dx.
\ee
Consequently, the integrals are in involution, i.e.,
\be
\{ \Psi_n, \Psi_m \} = 0,
\ee
if and only if there exists a function
$S_{nm} (u)$ such that
\be
{\pa \psi_n \over \pa u^i}
\eta^{ij} {\pa^2 \psi_m \over \pa u^j \pa u^k} =
{\pa S_{nm} \over \pa u^k},
\ee
i.e., if and only if the consistency relation
\be
{\pa \over \pa u^l} \left ( {\pa \psi_n \over \pa u^i}
\eta^{ij} {\pa^2 \psi_m \over \pa u^j \pa u^k} \right ) =
{\pa \over \pa u^k} \left ( {\pa \psi_n \over \pa u^i}
\eta^{ij} {\pa^2 \psi_m \over \pa u^j \pa u^l} \right ) \label{sovm2}
\ee
is fulfilled.
The consistency relation (\ref{sovm2}) is equivalent to
the equations (\ref{nonl4}).

It is proved similarly that
the equations (\ref{locint}) are equivalent to the condition
\be
\{ \Psi_n, H \} = 0, \ \ \ \ H = \int h (u(x)) dx.
\ee

We note that the equations (\ref{nonl4n}) are
equivalent to the condition that $L$
integrals are in involution with respect to
an arbitrary Dubrovin--Novikov bracket (a nondegenerate
local Poisson bracket of hydrodynamic type).

\bc
The hamiltonian system (\ref{arb}), (\ref{arb2})
is local if and only if it is generated by a family
of $L + 1$ integrals in involution with respect to
the Poisson bracket given by the constant Hamiltonian operator
$M^{ij}_2$ (\ref{nonlme}):
\be
\Psi_n = \int \psi_n (u(x)) dx, \ \
1 \leq n \leq L, \ \ H = \int h (u(x)) dx, \ \
\{ \Psi_n, \Psi_m \} = 0, \ \ \{ \Psi_n, H \} = 0.
\ee
In this case the system (\ref{arb}), (\ref{arb2}) is
an integrable bi-Hamiltonian system.
\ec

The important special class of integrals in involution
is generated by the equations of associativity
of two-dimensional topological quantum field theory.

\bl
A function $\Phi (u^1, ..., u^N)$ generates a family
of $N$ integrals in involution with respect to
the Poisson bracket given by the constant Hamiltonian operator
$M^{ij}_2$ (\ref{nonlme}):
\be
I_n = \int {\pa \Phi \over \pa u^n} (u(x)) dx, \ \
\{ I_n, I_m \} = 0, \ \  1 \leq n, m \leq N,
\ee
if and only if the function $\Phi (u)$ is a solution
of the equations of associativity in
two-dimensional topological quantum field theory
\be
\sum_{m=1}^N \sum_{n=1}^N \eta^{mn} {\pa^3 \Phi \over
\pa u^i \pa u^j \pa u^m} {\pa^3 \Phi \over \pa u^n \pa u^k \pa u^l}
=
\sum_{m=1}^N \sum_{n=1}^N \eta^{mn} {\pa^3 \Phi \over
\pa u^i \pa u^k \pa u^m}
{\pa^3 \Phi \over \pa u^n \pa u^j \pa u^l}. \label{nonl6}
\ee
\el

\section{Flat submanifolds with flat normal bundle, \\
Hessians and nonlocal Poisson brackets} \label{poverkh}

The nonlinear equations (\ref{nonl4}) and (\ref{nonl5})
describing all nonlocal Hamiltonian operators of
hydrodynamic type with flat metrics are exactly equivalent
to the conditions that a flat $N$-dimensional submanifold
with flat normal bundle, with the first fundamental form
$\eta_{ij} d u^i d u^j$ and the second fundamental
forms $\omega_n (u)$ given by Hessians of $L$ functions
$\psi_n (u),$ $1 \leq n \leq L$:
$$
\omega_n (u) = {\pa^2 \psi_n \over \pa u^i \pa u^j} d u^i d u^j,
$$
is embedded in an $(N + L)$-dimensional pseudo-Euclidean space.
In particular, the equations (\ref{nonl5}) are the Gauss
equations for this case, and the equations (\ref{nonl4})
are the Ricci equations
(the Peterson--Codazzi--Mainardi equations are
fulfilled automatically in this case).
By the Peterson--Bonnet theorem, any solution
$\psi_n (u),$ $1 \leq n \leq L,$ of the nonlinear system of equations
(\ref{nonl4}) and (\ref{nonl5}) defines a unique
up to motions $N$-dimensional submanifold
with flat normal bundle, with the first fundamental form
$\eta_{ij} d u^i d u^j$ and the second fundamental forms $\omega_n (u)$
given by Hessians of the functions $\psi_n (u)$, in
the $(N + L)$-dimensional pseudo-Euclidean space.

Let us consider an arbitrary $N$-dimensional flat submanifold
with flat normal bundle in an $(N + L)$-dimensional
pseudo-Euclidean space. Let $g_{ij} (u) d u^i d u^j$ be
the first fundamental form of this flat submanifold.
Then the following corollary of the theory
of submanifolds and the Bonnet theorem is true.

\bl
There locally exist functions $\Phi_n (u),$ $1 \leq n \leq L,$
such that the second fundamental forms of the submanifold under
consideration have the form
\be
\Omega_n (u) = \nabla_i \nabla_j \Phi_n d u^i d u^j,  \label{kvform}
\ee
where $\nabla_i$ is the operator of covariant
differentiation defined by the Levi-Civita connection of the
metric $g_{ij} (u)$. In this case the functions
$\Phi_n (u),$ $1 \leq n \leq L,$
satisfy the Ricci equations (\ref{nonl4n}) and the Gauss
equations (\ref{nonl5n}) for submanifolds
(the Peterson--Codazzi--Mainardi equations
are fulfilled automatically in this case).
Any solution of the system of equations (\ref{nonl4n})
and (\ref{nonl5n}) defines a unique up to motions
$N$-dimensional submanifold with flat normal bundle,
with the first fundamental form
$g_{ij} (u) d u^i d u^j$ and the second fundamental forms
$\Omega_n (u)$ (\ref{kvform}), in the $(N + L)$-dimensional
pseudo-Euclidean space.
\el

\section{Equations of associativity and nonlocal \\
Poisson brackets} \label{sass}

Here, we show that the class of nonlocal
Hamiltonian operators of hydrodynamic type
with flat metrics is nontrivial and describe
a rich and important family of operators of this class.
This family is generated by the equations of
associativity of two-dimensional topological quantum field theory.

Although the relations (\ref{nonl4}) and (\ref{nonl5})
differ essentially, they are very similar, and the case of
the natural reduction
when the relations
 (\ref{nonl4}) and (\ref{nonl5}) simply coincide is of particular
interest. Such reduction immediately leads to
the equations of
associativity of two-dimensional topological field theory.

In fact, the relations (\ref{nonl4}) and (\ref{nonl5})
coincide if
$L = N$, $\mu^{mn} = \eta^{mn}$
(or, for example, $\mu^{mn} = c \eta^{mn},$ where $c$
is an arbitrary nonzero constant) and
there exists a function $\Phi (u)$ such that
$\psi_n = {\pa \Phi / \pa u^n}$ for all $n$.
In this case both the relations (\ref{nonl4}) and (\ref{nonl5})
coincide with the equations of
associativity of two-dimensional topological field theory
for the potential $\Phi (u)$
(\ref{nonl6}).

Thus, any solution $\Phi (u)$ of the equations of
associativity (\ref{nonl6}), which are, as is well known,
consistent, integrable by the method of the inverse scattering problem
and possess a rich set of nontrivial solutions
(see \cite{3}), defines a nonlocal Hamiltonian operator
of hydrodynamic type with a flat metric:
\be
L^{ij} = \eta^{ij} {d \over dx} +
\sum_{m=1}^N \sum_{n =1}^N
\eta^{mn} \eta^{ip} \eta^{jr}
{\pa^3 \Phi \over \pa u^p \pa u^m \pa u^k}
u^k_x \left ( {d \over dx} \right )^{-1} \circ
{\pa^3 \Phi \over \pa u^r \pa u^n \pa u^s} u^s_x, \label{nonl7}
\ee
and, moreover, defines a pencil of compatible
Hamiltonian operators:
\be
L^{ij}_{\lambda_1, \lambda_2} = \lambda_1
\eta^{ij} {d \over dx} +
\lambda_2 \sum_{m=1}^N \sum_{n =1}^N
\eta^{mn} \eta^{ip} \eta^{jr}
{\pa^3 \Phi \over \pa u^p \pa u^m \pa u^k}
u^k_x \left ( {d \over dx} \right )^{-1} \circ
{\pa^3 \Phi \over \pa u^r \pa u^n \pa u^s} u^s_x, \label{nonl7a}
\ee
where $\lambda_1$ and $\lambda_2$ are arbitrary constants.
In particular, if the function $\Phi (u)$ is an arbitrary
solution of the equations of associativity
(\ref{nonl6}), then the operator
\be
L^{ij}_{0, 1} = \sum_{m=1}^N \sum_{n =1}^N
\eta^{mn} \eta^{ip} \eta^{jr}
{\pa^3 \Phi \over \pa u^p \pa u^m \pa u^k}
u^k_x \left ( {d \over dx} \right )^{-1} \circ
{\pa^3 \Phi \over \pa u^r \pa u^n \pa u^s} u^s_x \label{nonl7b}
\ee
is a Hamiltonian operator compatible with the
constant Hamiltonian operator
\be
L^{ij}_{1, 0} =
\eta^{ij} {d \over dx}.  \label{nonl7c}
\ee
Therefore, for any solution of the equations of associativity
(\ref{nonl6}) we obtain the corresponding integrable hierarchies
(see Section \ref{hier5}).

The metric
$\eta^{ij}$ always defines a nondegenerate
invariant symmetric bilinear form in
the corresponding Frobenius algebras, namely,
commutative associative algebras $A (u)$ in an
$N$-dimensional vector space with a basis $e_1, \ldots, e_N$
and the multiplication (see \cite{3})
\be
e_i \ast e_j = \eta^{ks}
{\pa^3 \Phi \over \pa u^s \pa u^i \pa u^j} e_k, \label{alg1}
\ee
\be
<e_i, e_j> = \eta_{ij}, \ \ \ \
<e_i \ast e_j, e_k> = <e_i, e_j \ast e_k>.  \label{alg2}
\ee
In this case the condition of associativity
\be
(e_i \ast e_j) \ast e_k =
e_i \ast (e_j \ast e_k)       \label{alg3}
\ee
in the algebras $A (u)$ is equivalent to the equations (\ref{nonl6}).
We recall (see Dubrovin, \cite{3}) that locally
on the tangent space in every point of any Frobenius manifold
there is a structure of Frobenius algebra
(\ref{alg1}), (\ref{alg2}),
(\ref{alg3}) defined by a certain solution
of the equations of associativity
(\ref{nonl6}) and depending on the point smoothly (besides,
it is also required the fulfilment of additional conditions on
Frobenius manifolds but we do not consider these conditions here;
roughly speaking, it is required the presence of a unit and
the quasihomogeneity).
Thus, for every Frobenius manifold
there are nonlocal Hamiltonian operators of the form (\ref{nonl7})
and pencils
(\ref{nonl7a}) connected to the manifold.

We have considered the nonlocal Hamiltonian operators of
the form (\ref{nonlm}) with flat metrics and came to
the equations of associativity defining the affinors of such operators.
A statement that is in some sense the converse is also true,
namely, if all the affinors $w_n$
of a nonlocal Hamiltonian operator (\ref{nonlm})
with $L = N$ are defined by an arbitrary solution
$\Phi (u)$ of the equations of associativity
\be
\sum_{m=1}^N \sum_{n=1}^N \mu^{mn} {\pa^3 \Phi \over
\pa u^i \pa u^j \pa u^m} {\pa^3 \Phi \over \pa u^n \pa u^k \pa u^l}
=
\sum_{m=1}^N \sum_{n=1}^N \mu^{mn} {\pa^3 \Phi \over
\pa u^i \pa u^k \pa u^m}
{\pa^3 \Phi \over \pa u^n \pa u^j \pa u^l}
\ee
by the formula
$$
(w_n)^i_j (u) = \zeta^{is} \xi^r_j {\pa^3 \Phi \over
\pa u^n \pa u^s \pa u^r},
$$
where $\zeta^{is},$ $\xi^r_j$ are arbitrary nondegenerate
constant matrices, then the metric of this Hamiltonian
operator must be flat. But, in general, it is not necessarily that
this metric will be constant in the local coordinates under consideration.

The structural flows (see \cite{1}, \cite{4})
of the nonlocal Hamiltonian operator (\ref{nonl7})
have the form:
\be
u^i_{t_n} = \eta^{is} {\pa^3 \Phi \over \pa u^s \pa u^n \pa u^k}
u^k_x.  \label{nonl8}
\ee
These systems are integrable bi-Hamiltonian systems of
hydrodynamic type and coincide with the primary part of
the Dubrovin hierarchy constructed by any solution
of the equations of associativity in \cite{3}.
The condition of commutation for
the structural flows (\ref{nonl8})
is also equivalent to the equations of associativity (\ref{nonl6}).

From the consideration of the previous Section
we have
that by the Peterson--Bonnet theorem any solution
$\Phi (u)$ of the equations of associativity (\ref{nonl6})
defines a unique up to motions
$N$-dimensional submanifold
with flat normal bundle, with the first fundamental form
$\eta_{ij} d u^i d u^j$ and the second fundamental forms
$$
\omega_n (u) = {\pa^3 \Phi \over \pa u^n \pa u^i \pa u^j} d u^i d u^j
$$
given by the potential $\Phi (u)$, in a $2N$-dimensional
pseudo-Euclidean space. Here the signature of the
ambient pseudo-Euclidean space is completely defined
by the signature of the metric $\eta_{ij}$.
In this case the equation of associativity coincides with
both the Gauss equation and the Ricci equation of
the embedded submanifold.
Thus, a special class of flat
$N$-dimensional
submanifolds with flat normal bundle in a $2N$-dimensional
pseudo-Euclidean space locally corresponds to
$N$-dimensional Frobenius manifolds.

A great number of concrete examples of Frobenius
manifolds and solutions of the equations of associativity
is given in Dubrovin's paper \cite{3}.
Here we do not give a great number of examples and
an explicit form of the corresponding Hamiltonian
operators, flat submanifolds and integrable hierarchies
to save place and consider only one simple example from
\cite{3} as an illustration.
Let $N = 3$ and the metric $\eta_{ij}$ be antidiagonal
\be
(\eta_{ij}) =
\left ( \begin{array} {ccc} 0&0&1\\
0&1&0\\
1&0&0
\end{array} \right ),
\ee
and the function $\Phi (u)$ has the form
$$
\Phi (u) = {1 \over 2} (u^1)^2 u^3 +
{1 \over 2} u^1 (u^2)^2 + f (u^2, u^3).
$$
In this case $e_1$ is the unit in the Frobenius algebra
(\ref{alg1}), (\ref{alg2}),
(\ref{alg3}), and the equation of associativity
(\ref{nonl6}) for the function $\Phi (u)$ is equivalent to
a remarkable integrable equation of Dubrovin for the function
$f (u^2, u^3)$:
\be
{\pa^3 f \over \pa (u^3)^3} = \left (
{\pa^3 f \over \pa (u^2)^2 \pa u^3} \right )^2 -
{\pa^3 f \over \pa (u^2)^3} {\pa^3 f \over \pa u^2 \pa (u^3)^2}.
\label{f}
\ee
This equation is connected to quantum cohomology of
projective plane and classical problems of
enumerative geometry (see \cite{6}). In particular,
all nontrivial polynomial solutions of the equation (\ref{f})
that satisfy the requirement of the quasihomogeneity
and, consequently, locally define a structure of
Frobenius manifold are described by Dubrovin in \cite{3}:
\be
f = {1 \over 4} (u^2)^2 (u^3)^2 + {1 \over 60} (u^3)^5,\ \ \
f = {1 \over 6} (u^2)^3 u^3 + {1 \over 6} (u^2)^2 (u^3)^3 +
{1 \over 210} (u^3)^7,
\ee
\be
f = {1 \over 6} (u^2)^3 (u^3)^2 + {1 \over 20} (u^2)^2 (u^3)^5 +
{1 \over 3960} (u^3)^{11}.
\ee

As is shown by the author in \cite{7} (see also \cite{8}),
the equation (\ref{f}) is equivalent to the integrable
nondiagonalizable homogeneous system of hydrodynamic type
\be
\left ( \begin{array} {c} a\\ b\\ c
\end{array} \right )_{u^3} =
\left ( \begin{array} {ccc} 0 & 1 & 0\\  0 & 0 & 1\\
- c & 2b & - a
\end{array} \right )  \left ( \begin{array} {c} a\\ b\\ c
\end{array} \right )_{u^2},
\ee
\be
a = {\pa^3 f \over \pa (u^2)^3}, \ \ \
b = {\pa^3 f \over \pa (u^2)^2 \pa u^3},\ \ \
c = {\pa^3 f \over \pa u^2 \pa (u^3)^2}.
\ee
In this case the affinors of the nonlocal Hamiltonian
operator
(\ref{nonl7}) have the form:
\be
(w_1)^i_j (u) = \delta^i_j, \ \ \
(w_2)^i_j (u) =
\left ( \begin{array} {ccc} 0 & b & c\\  1 & a & b\\
0 & 1 & 0
\end{array} \right ),\ \ \
(w_3)^i_j (u) = \left ( \begin{array} {ccc}
0 & c & b^2 - a c\\  0 & b & c\\
1 & 0 & 0
\end{array} \right ).
\ee

\medskip

\begin{flushleft}
Center for Nonlinear Studies,\\
L. D. Landau Institute for Theoretical Physics,\\
Russian Academy of Sciences\\
e-mail: mokhov@mi.ras.ru; mokhov@landau.ac.ru; mokhov@bk.ru\\
\end{flushleft}

\end{document}